\documentclass[a4paper,12pt]{report}
 \usepackage[english]{babel}
            \usepackage{graphicx}
\usepackage{hyperref}
\usepackage{parskip}

\usepackage{varioref}
\usepackage{amsmath}
\usepackage{amsfonts}
\usepackage{amssymb}
\usepackage{amsthm}
\usepackage{color}
\usepackage{makeidx}
\usepackage{xspace}
\usepackage{tabularx}
\usepackage{verbatim}

\usepackage{geometry}
\usepackage{texnames}
\usepackage{cite}

\begin{document}

\newcommand{\opL}{{\mathbf L}}
\newcommand{\opLo}{{\mathbf L}_1}
\newcommand{\opB}{\mathbf B}
\newcommand{\opA}{\mathbf A}
\newcommand{\opBo}{{\mathbf B}_1}
\newcommand{\opBt}{{\mathbf B}_2}
\newcommand{\opAo}{{\mathbf A}_1}
\newcommand{\opAt}{{\mathbf A}_2}

\newcommand{\beq}{\begin{equation}}
\newcommand{\eeq}{\end{equation}}

MSC 47A50, 47A52 \\
\begin{center}
{\bf \Large Nonclassic boundary value problems in the theory of
irregular systems of equations with partial derivatives\footnote{This work is fuilfilled
within International Science and Technology
Cooperation Program of China and Russia, project 2015DFA70850,  NSFC Grant No. 61673398.
 }}
$$\,$$
{\large Nikolai~{Sidorov and}} 
{\large Denis~{Sidorov}}\\
\end{center}

 $$\,$$
{\bf Abstract}. {\it The linear PDE
$\,\,\opB \opL (\frac{\partial}{\partial x}) u  =\opL_1(\frac{\partial}{\partial x})u +f(x)$
with nonclassic conditions on boundary $\partial \Omega$
is considered. Here
 $\opB$ is linear noninvertible bounded operator acting from linear space $E$ into $E,$
$x=(t,x_1,\dots, x_m) \in \Omega, $ $\Omega \subset {\mathbb R}^{m+1}.$
It is assumed that $\opB$
 enjoys the skeleton decomposition
$\opB={\mathbf A}_1 {\mathbf A}_2,$ ${\mathbf A}_2 \in {\mathcal L}(E\rightarrow E_1),$  ${\mathbf A}_1 \in {\mathcal L}(E_1\rightarrow E)$ where
$E_1$ is linear normed space. Differential operators
$\opL, \, \opL_1$ are as follows
$${\mathbf L}\left (\frac{\partial}{\partial x}\right) = \frac{\partial^n}{\partial t^n} + \sum\limits_{k_0+k_1+ \dots+ k_m \leq n-1} a_{k_0 \dots k_m}(x) \frac{\partial^{k_0+\dots + k_m}}{\partial t^{k_0} \partial x_1^{k_1} \dots x_m^{k_m} }, $$ 
$${\mathbf L}_1\left (\frac{\partial}{\partial x}\right) =  \sum\limits_{k_0+k_1+ \dots+ k_m \leq n_1} b_{k_0 \dots k_m}(x) \frac{\partial^{k_0+\dots + k_m}}{\partial t^{k_0} \partial x_1^{k_1} \dots x_m^{k_m} },\, n_1<n $$
where coefficients $ a_{k_0 \dots k_m}: \Omega \rightarrow {\mathbb R}^1,$ $b_{k_0 \dots k_m}: \Omega 
\rightarrow {\mathbb R}^1 $ are sufficiently smooth and defined in  $\overline{\Omega}, 0\in \Omega.$
In the concrete cases the domains of definition of operators $\opL, \opLo$  consist of 
linear manifolds $E_{\partial}$ of sufficiently smooth abstract functions $u(x)$ with domain in $\Omega$ and their ranges in  $E,$ which satisfy certain system of homogeneous boundary conditions. 
 The abstract function  $f: \Omega \subset {\mathbb R}^{m+1} \rightarrow E $ of
argument $(t,x_1, \dots, x_m)$
 is assumed to be given.  It is requested to find the solution $u: \Omega \subset {\mathbb R}^{m+1} \rightarrow E_{\partial},$ abstract function $u(x)$ satisfy certain condition on boundary $\partial \Omega.$ 
The concept of a skeleton chains is introduced as sequence of 
linear operators $\opB_i \in {\mathcal L}(E_i \rightarrow E_i), \, i=1,2,\dots, p,$
where 
$E_i$ are
 linear spaces corresponding to the skeleton decomposition of operator $\opB.$
It is assumed that irreversible operator  $\opB$
generates skeleton chain of the finite length $p.$ The problem is reduced to a regular split system
with respect to higher order derivative terms with certain initial and boundary conditions.
}






\section*{Introduction}

Let linear bounded operator ${\mathbf B}$ acting from linear space 
$E$ to $E$ has no  inverse operator.  Differential operators
$${\mathbf L}\left (\frac{\partial}{\partial x}\right) = \frac{\partial^n}{\partial t^n} + \sum\limits_{k_0+k_1+ \dots+ k_m \leq n-1} a_{k_0 \dots k_m}(x) \frac{\partial^{k_0+\dots + k_m}}{\partial t^{k_0} \partial x_1^{k_1} \dots x_m^{k_m} }, $$ 
$${\mathbf L}_1\left (\frac{\partial}{\partial x}\right) =  \sum\limits_{k_0+k_1+ \dots+ k_m \leq n_1} b_{k_0 \dots k_m}(x) \frac{\partial^{k_0+\dots + k_m}}{\partial t^{k_0} \partial x_1^{k_1} \dots x_m^{k_m} },\, n_1<n $$
are defined. Here  coefficients $ a_{k_0 \dots k_m}: \Omega \subset {\mathbb R}^{m+1}\rightarrow {\mathbb R}^1,$ $b_{k_0 \dots k_m}: \Omega \subset {\mathbb R}^{m+1}
\rightarrow {\mathbb R}^1$ are sufficiently smooth and defined in $\overline{\Omega}, 0\in \Omega.$
The domains of definition of operators $\opL, \opLo$ 
 consist of 
linear manifolds $E_{\partial}$ sufficiently smooth functions in $\Omega$ with their ranges in  $E,$ which satisfy certain system of homogeneous boundary conditions. 

Abstract function $f: \Omega \subset {\mathbb R}^{m+1} \rightarrow E $
of argument $x$ is assumed to be given and the problem is to find the solution $u: \Omega \subset {\mathbb R}^{m+1} \rightarrow E_{\partial}$
which satisfy  linear PDE
\begin{equation}
\opB \, \opL \left ( \frac{\partial}{\partial x}\right ) u = \opLo \left ( \frac{\partial}{\partial x}\right ) u +f(x). \label{eq1}
\end{equation} 
Operator $\opB$ is assumed independent of $x.$
If operator $\opB$ has inverse bounded operator then Eq.~\eqref{eq1} is called {\it regular} and
otherwise it is called {\it irregular} equation.  If $E = {\mathbb R}^N$
and $\det B \neq 0,$ then Eq.~\eqref{eq1}  is the system of linear partial differential equations (PDE) of Kovalevskaya type, and we have a well known regular problem of the PDE theory. The foundation of many branches of modern general theory of PDE systems was constructed
by I.G.~Petrovskii~\cite{petro}.
In regular case  the initial conditions for Eq. \eqref{eq1} can be defined as follows
\begin{equation}
\frac{\partial^i u}{\partial t^i} \bigg|_{t=0} = \varphi_i(x_1, \dots , x_m), \, i=0,1, \dots, n-1.
 \label{eq2}
\end{equation}
Here functions $\varphi_i$ are analytical functions in $\Omega.$  
If $f$ is analytic function in $t,$  $x_1, \dots , x_m$ in $\Omega,$
then Cauchy problem \eqref{eq1} -- \eqref{eq2} is not only solvable but  also {\it well-posed}
in class of analytic functions.

The well-posedness of Cauchy problem is challenging issue  even for linear 
PDE systems in  spaces of non-analytic functions.  They are usually solved in  class of functions satisfying certain estimates~\cite{petro}.
As result of publication of  S.L.~Sobolev work~\cite{sobolev} where he introduced 
nowadays called Sobolev equations, I.G.~Petrovsky and L.A.~Lusternik on their seminars
have attracted mathematicians attention to PDE systems unresolved with respect to the highest time derivative.

It is to be noted that there are numerous new models in the natural sciences, engineering, and mathematical economics (here readers may refer to \cite{appl, appl2, appl3} and other) formulated in terms of  systems of Eqs~\ref{eq1}.

Irregular models enable study of systems behavior in critical situations. 
At present, the basis of relevant theory is constructed for certain classes of equations. For example,  the theory and numerical methods for differential-algebraic equations has
been constructed.

  The intensive studies of more complex theory of irregular PDE and abstract irregular differential operator
equations are coducted but there are still a lot of unexplored
problems.

 If $\opB$ is normally solvable operator, $x \in {\mathbb R}^1$ then the approach in the theory of Eq. \eqref{eq1} can be 
 based on the expansion of the Banach space into direct sum in accordance with  Jordan structure of operator $\opB$~\cite{sid_mono}-\cite{main_mono}  and some results
from the theory of semigroup with kernels\cite{svir}. These approaches are already 
employed for various problems of modern mathematical modeling, see e.g.~\cite{appl, appl2,appl3}.

In this field the analytical methods were proposed for constructing classical and generalized solutions of  Cauchy problem for ordinary operator-differential equations for $ x \in {\mathbb R}^1 $ in Banach spaces with irreversible operator in the main part.

The theory of irregular operator-differential PDEs in Banach spaces in the multi-dimensional case for $ x \in {\mathbb R}^n, \, n \geq 2 $ to be constructed. 
There are only initial results  in this field published in preprints~\cite{sid_blag, sid_rom_blag}. Therefore, the construction of the general theory of Eq.~\eqref{eq1} with irreversible operator $\opB$  is of theoretical interest. It is also important for the modern mathematical models based on 
irregular systems.

It is to be outlined that classical initial
Cauchy conditions  \eqref{eq2} for Eq.~\eqref{eq1} play very limited role. 
Indeed, because of the irreversibility
the operator $\opB$  time direction is characteristic and
 functions $ \varphi_i $ can not be  arbitrary selected in the initial conditions \eqref{eq2}! 
Then appear 
the question of reasonable formulation and methods of solution of non-classic boundary
problems for the system \eqref{eq1}, taking into account the structure of the operator $ \opB. $
The objective of present work is to solve this problem.
The similar Goursat problem was addressed in \cite{sid_blag91} using other
methods. In sections 2 and 3 of present work this problem is solved for
irreversible operator $ \opB $ which enjoy skeleton decomposition
$\opB = \opAo \opAt,$ где $\opAo \in {\mathcal L}(E\rightarrow E_1),$
$\opAt \in {\mathcal L}(E_1\rightarrow E),$ where
$E_1$ is normed space.

The remainder of the paper is organized as follows. 
Sec.~0.1 presents introductory  concerning the skeleton chains of linear operators 
using results\cite{arxiv}. The concept of a regular and singular skeleton chains
is introduced. It is proved that operator $\opB$ must be nilpotent  in case of singular skeleton chain. In Sec.~0.2 it is assumed that noninvertible operator $\opB$
 generates a skeleton chain of linear operators of finite length
 $p$ and it is demonstrated that irregular Eq.~\eqref{eq1} can be reduced to
the recurrent sequence of  $p+1$ equations. It is to be noted that each
equation of this sequence is regular under the
natural restrictions on differential operators $\opL, \opLo$
and  certain initial-boundary conditions. Therefore if operator
$\opB$ has skeleton chain  of length $p$ then solution of irregular 
Eq.~\eqref{eq1}  can be reduced to regular system from  $p+1$-th equation.

Proposed approach can be employed for wide range of concrete problems~\eqref{eq1} due to finite 
length of skeleton chain of  finite-dimensional operator $\opB.$

The formulas relating the solution of Eq.~\eqref{eq1} with the solution of
reduced regular system are derived.  This result allows us in Sec.~0.2  to set new well-posed non-classic boundary conditions for Eq. \eqref{eq1} for which the equation enjoy unique solution as demonstrated in Sec.~0.3 -- 0.5. 
For applications, it is important that this solution can be found by solving the sequence of regular problem proposed in this paper. Corresponding results and examples are given in Sec.~0.3 and Sec.~0.4. Finally Sec.~0.5 demonstrates the results generalization to  non-linear equations.

 \section{Skeleton chains of linear operator}   

 Let $\opB \in  {\mathcal L}(E \rightarrow E), $ and $\opB = \opAo \opAt, $
where $\opAt \in {\mathcal L}(E \rightarrow E_1), \opAo \in {\mathcal L}(E_1 \rightarrow E),$
$E_1, E$ are linear normed spaces. The following definitions
can be introduced.

Decomposition $\opB = \opAo \opAt$
is called  {\it skeleton decomposition of operator} $\opB.$  
Let us introduce linear operator $\opBo = \opAt \opAo.$
Obviously $\opBo \in {\mathcal L}(E_1 \rightarrow E_1).$  
If operator $\opBo$ has bounded inverse or it is null operator acting from $E_1$ to $E_1,$ 
then  $\opB$  generate
{\it skeleton chain} $\{\opBo \}$ of length 1. Then operator $\opBo$ can be called as
{\it skeleton-attached operator} to operator $\opB.$
This chain is called {\it singular } if  $\opBo = 0$  and 
regular if  
$\opBo \neq 0.$ If $\opBo$ is 
irreversible non-null operator then it is assumed to have skeleton decomosition
$\opBo = {\mathbf A}_3 {\mathbf A}_4,$ где $  {\mathbf A}_4 \in  {\mathcal L}(E_1 \rightarrow E_2), $ $  {\mathbf A}_3 \in  {\mathcal L}(E_2 \rightarrow E_1), $  where
$E_2$ is new linear normed space. Obviously in this case $\opAt \opAo  =  {\mathbf A}_3 {\mathbf A}_4$  and operator $\opBt =  {\mathbf A}_4 {\mathbf A}_3 \in  {\mathcal L}(E_2 \rightarrow E_2)$ can be introduced.
If it turns out that $\opBt$ has bounded inverse or  $\opBt \equiv 0,$ then
 $\opB$ has skeleton chain $\{ \opBo, \opBt \}$ of 
length 2.
Chain $\{\opBo, \opBt \}$ is generate if $\opBt=0$
 and nongenerate otherwise.

Therefore this process can be continued for a number of linear operators
classes by introduction of the normed linear spaces $E_i,\, i=1,\dots, p$  
and by bounded operators construction
$  {\mathbf A}_{2i} \in  {\mathcal L}(E_{i-1} \rightarrow E_i),$
$  {\mathbf A}_{2i-1} \in  {\mathcal L}(E_{i} \rightarrow E_{i-1}),$
which satisfy the following equalities
\begin{equation}
{\mathbf A}_{2i} {\mathbf A}_{2i-1} = {\mathbf A}_{2i+1} {\mathbf A}_{2i+2},
\, i=1,2,\dots, p-1.
\label{eq3}
\end{equation}
 Eq. \eqref{eq3} defines sequence of linear operators $\{ {\mathbf B}_{1}, \dots,  {\mathbf B}_{p}\}$ as follows
\beq
{\mathbf B}_{i} = {\mathbf A}_{2i} {\mathbf A}_{2i-1}, \, i=1,2, \dots , p.
\label{eq4}
\eeq
Obviously $ {\mathbf B}_{i} \in {\mathcal L}(E_i\rightarrow E_i).$
Here  operator  ${\mathbf B}_{p}$ either has inverse bounded or 
${\mathbf B}_{p}$  is null operator acting from $E_p$ to $E_p.$
This process can be formalized as following definition.\\

{\it Definition 1.}\\
Let $\opB = \opAo \opAt$ and operators $\{{\mathbf A}_{i}\}_{i=1}^{2p}$
satisfy equality~\eqref{eq3}.  Let operators 
 $\{ \opBo, \dots, {\mathbf B}_{p} \}$ defined by formula~\eqref{eq4},
 and operators   $\{ \opBo, \dots, {\mathbf B}_{p-1} \}$ are noninvertible, 
and operator   ${\mathbf B}_{p}$  has bounded inverse or null operators acting from $E_p$ to $E_p.$
Then operator $\opB$ {\it generates skeleton chain} of linear operators
 $\{ \opBo, \dots, {\mathbf B}_{p} \}$ on lenght $p.$
If 
 ${\mathbf B}_{p}  \neq 0$  then the chain is {\it regular}, 
if ${\mathbf B}_{p}  = 0$ then the chain is called {\it singular.} 
Operators   
$\{ \opBo, \dots, {\mathbf B}_{p} \}$ called skeleton-attached to operator $\opB.$

 The most important linear opeators generating skeleton chains of the finite lenghts are considered below:

 \begin{enumerate}
\item Let $E= {\mathbb R}^n,$  then square matrix $\opB: {\mathbb R}^m \rightarrow {\mathbb R}^m$  with $\det \opB =0$ obviously has skeleton chain  $\{\opB_1, \dots , \opB_p \}$ 
of decreasing dimentions. The final matrix $\opB_p$ will be 
regular or null matrix, $\det \opB_i = 0, \, i=1,\dots,p-1.$

\item Let $E$ is infinite dimentional normed space, then
finite operator  $\opB = \sum\limits_{i=1}^n \langle \cdot, \gamma_i \rangle z_i,$  where $\{z_i \} \in E,$ $\gamma_i \in E^*$ has skeleton chain consisting from finite number  of matrices $\{ \opBo, \dots, {\mathbf B}_{p} \} $ of decreasing dimentions. Here $\opB_1 = || \langle z_i, \gamma_j \rangle ||_{i,j=1}^n$ is first element of this chain,
$\det \opB_i=0, \, i=1,\dots, p-1.$
$\opB_p$ is null matrix or $\det B_p \neq 0.$
\end{enumerate}
Here according definition 1 length of the chain $p=1$  if
$ \det  [\langle z_i, \gamma_j \rangle]_{i,j=1}^n \neq 0$ or  $\langle z_i, \gamma_j \rangle = 0, i,j=1,2, \dots, n.$  In general case the chain 
is always 
consists from finite number of matrices.

Using \eqref{eq3}, \eqref{eq4}  and Definition 1 the following result can be formulated.\\

\noindent {\bf Lemma 1.}\\
If operator  $\opB$ has skeleton chain of lenght  $p$ then 
\beq
\opB^n = \opAo {\mathbf A}_3 \dots  {\mathbf A}_{2n-1}  {\mathbf B}_{n-1}
 {\mathbf A}_{2n-2} {\mathbf A}_{2n-4} \dots \opAt,\,\, n=1,\dots,p+1,
\label{eq5}
\eeq
where $\opBo, \opBt, \dots {\mathbf B}_p $ are elements of the skeleton chain of operator $\opB.$\\

From Lemma 1 it follows\\
$\,$\\
{\bf Corollary.} \\
If operator  $\opB$ has singular skeleton chain of length $p,$  then $\opB$  is nilpotent 
operator of index  $p+1.$\\

To proof the Corollary it is sufficient to put $n=p+1$ and demonstrate that degree
 $\opB^{p+1}$ is null operator because  $\opB_p$ is null operator 
due to above introduced definition of singular  skeleton chain.\\

\section{Reduction of abstract irregular Eq.~\eqref{eq1} to\\ the sequence of regular equations}

We always assume that operator $\opB$ and  linear operators $\{ \opA_i\}_{i=1}^{2p}$ from skeleton chain of operator  $\opB$ be independent on $x$ and
commutative with linear operators $\opL$ and $\opL_1.$
In this paragraph for sake of clarity it is assumed that operators $\opL$ and $\opL_1$
can be different from above introduced differential  operators $\opL(\frac{\partial }{\partial x}),$
$\opL_1(\frac{\partial }{\partial x})$ and 
equation can be considered in abstract form
\beq
\opB \opL u = \opL_1 u +f.
\label{eq7}
\eeq
Eq.~\eqref{eq1}  can be considered as special case of Eq.~\eqref{eq7}. Obviously, the introduced commutativity 
condition is fuilfilled for Eq.~\eqref{eq1}  with linear operator  $\opB$
independent of $x$ and introduced differential operators $\opL(\frac{\partial }{\partial x}),$
$\opL_1(\frac{\partial }{\partial x}).$
 
Let us reduce Eq.~\eqref{eq7} to system of  $p+1$ equations which are regular 
in certain conditions imposed on operators  $\opL, \opL_1.$
Let us start with simple case   $p=1.$
We introduce system of two equation
\beq
\opB_1 \opL u_1 = \opL_1 u_1 +\opA_2 f,
\label{eq8}
\eeq
\beq
\opL_1 u = -f +\opA_1 \opL u_1.
\label{eq9}
\eeq
where  $u\in E$ and $u_1 \in E_1.$
The decomposed system \eqref{eq8} -- \eqref{eq9} can be obtained by formal multiplication of Eq.~\eqref{eq7}  on operator  $\opA_2$  from the skeleton decomposition of operator  $\opB$
and making notation $u_1 = \opA_2 u.$

It is to be noted that system  \eqref{eq8} -- \eqref{eq9}  is splited and  $\opB_1$ is
invertible operator. Therefore if operators  $\opB_1 \opL -  \opL_1,$ and operator $\opL_1$
have bounded inverse operators then the unique solution can be constructed.
Of course without additional conditions there remains an open question: 
{\it is constructed solution $u(x)$ satisfy  Eq.~\eqref{eq7}?}

Let us introduce two lemmas establishing the link between Eq.~\eqref{eq7}
 ans system \eqref{eq8} -- \eqref{eq9}.\\

{\bf Lemma 2.}\\
Let  $u^*$  satisfy Eq.~\eqref{eq7} and operator 
 $\opL_1$  has left inverse. Then pair $u_1^* = A_2 u^*, u^*$
satisfy system \eqref{eq8} -- \eqref{eq9}.

\begin{proof}
Based on conditions of the Lemma the following equality is satisfyed
\beq
\opAo \opAt \opL u^* = \opL_1 u^* +f.
\label{eq11}
\eeq
From \eqref{eq11} because of commutativity condition the following equality is valid
\beq
\opAo \opL \opA_2 u^* = \opL_1 u^* +f.
\label{eq12}
\eeq
 and
\beq
\opAt \opAo \opL \opA_2 u^* = \opL_1 \opA_2 u^* +\opA_2 f.
\label{eq13}
\eeq
The latter equality demonstrates that   $u_1^* = \opA_2 u^*$  is solution of 
Eq.~\eqref{eq8}.  Substitution $u_1^*$ into the right hand side of Eq.~\eqref{eq9} 
yields the following equation with known right hand side  $$\opL_1u= -f+ \opA_1 \opA_2 \opL u^*$$ 
with respect to $u.$ Here the operators commutativity property is employed. 
The solution exists for such equation. Indeed,  due to Eq.~\eqref{eq11}
right hand side of the equation is equal to  $\opL_1 u^*.$ 
Hence for given $u^*$ and $u_1^*=\opA_2 u^*$ the right hand side belongs
to the range of operator  $\opL_1.$
Therefore,  $\opL_1u = \opL_1 u^*.$  Because operator $\opL_1$ has left inverse,
then $u^*$  is unique solution ro Eq.~\eqref{eq9} for $u_1=\opA_2u^*$.  
Lemma 2 is proved.\\
\end{proof}

{\bf Lemma~3.}\\
Suppose that there exists a pair $(u_1^*, u^*)$ solution to Eq.~\eqref{eq8}, \eqref{eq9}.
Let operator $\opL_1$ has right inverse $\opL_1^{-1}.$
 Then element  $u^*$ satisfy Eq.~\eqref{eq7}.\\

\begin{proof}
Element $-f + \opA_1 \opL u_1^*$ belongs to the range of operator 
$\opL_1,$ because pair $u_1^*, u^*$ satisfies system of Eqs.~\eqref{eq8} -- \eqref{eq9}.
Hence $u^* = \opL_1^{-1} (-f + \opA_1 \opL u_1^*)$ because $-f + \opA_1 \opL u_1^* \in {\mathcal R}(\opL_1),$ where 
${\mathcal R}(\opL_1) $ is the range of operator $\opL_1.$
It is to be demonstrated that constructed element $u^*$ 
satisfy Eq.~\eqref{eq7} because by hypothesis
 $u_1^*$ satisfy Eq.~\eqref{eq8}.
Indeed, substitution of the constructed $u^*$ into Eq.~\eqref{eq7},
where $\opB = \opA_1 \opA_2,$
yields
$$\opA_1 \opA_2 \opL \opL_1^{-1} (-f+ \opA_1 \opL u_1^*) = \opA_1 \opL u_1^*. $$
Taking into account operators commutativity the following equality is valid
\beq
\opA_1 \opL \{ \opA_2 \opL_1^{-1} (-f+\opA_1 \opL u_1^*)-u_1^* \}=0
\label{eq14}
\eeq
 Because $\opA_1$ and $\opL$ are linear operators then
it remains to verify  in Eq.~\eqref{eq14}  element in braces is zero.
Since $u_1^* \in E_1$ satisfy Eq.~\eqref{eq8},
where $\opB_1 = \opA_2 \opA_1,$
then the following equality is valid
\beq
\opA_2 (\opA_1 \opL u_1^* - f ) = \opL_1 u_1^*.
\label{eq15}
\eeq
Hence $\opA_2 (\opA_1 \opL u_1^* - f) \in {\mathcal R}(\opL_1)$
and $u_1^* = \opL_1^{-1} \opA_2 (\opA_1 \opL u_1^* -f),$
where $\opL_1^{-1}$ is right inverse to operator  $\opL_1.$
 Because $\opA_2 \opL_1 = \opL_1 \opA_2$
then
$\opA_2 = \opL_1^{-1} \opA_2 \opL_1.$
This yields $\opA_2 \opL_1^{-1} = \opL_1^{-1} \opA_2 \opL_1 \opL_1^{-1}.$
From equality  $\opL_1 \opL_1^{-1} = {\mathbf I}$ it follows
that right inverse $\opL_1^{-1}$ also commutative with operator $\opA_2.$
 Then Eq.~\eqref{eq15} can be represented as
$$\opA_2 \opL_1^{-1} (\opA_1 \opL u_1^* -f )- u_1^*=0. $$
Thus we have shown that expression in braces in Eq.~ \eqref{eq14} is 
zero. Lemma 3 is proved.
\end{proof}

Let us concentrate on general case of skeleton chain of arbitrary finite lenght $p.$
We assume operator  $\opB$ to have skeleton chain
$\{ \opB_1, \dots , \opB_p  \},$ $ p \geq 1,$
and linear operators $\{A_i \}_{i=1}^{2p}$
represent decomposition of $\opB_i$ which skeleton attached to $\opB.$
Introduce
\beq
u_i = \prod_{j=1}^i \opA_{2j} u, \, i =1,\dots, p,
\label{eq16}
\eeq
where
$u_i \in E_i, \,u_i = \opA_{2i} u_{i-1}, u_0:=u.$

If $u_0$ satisfy Eq.~\eqref{eq7} then for
 $p\in {\mathbb N}$ by {\it Definition 1}  we get  equalities 
\beq
\opB_p \opL u_p = \opL_1 u_p +\prod_{j=1}^p \opA_{2j} f,
\label{eq16}
\eeq

\beq
\opL_1 u_i = -\prod_{j=1}^i \opA_{2j} f + \opA_{2i+1} \opL u_{i+1},
\label{eq17}
\eeq

\beq
\opL_1 u = -f + \opA_1 \opL u_1.
\label{eq18}
\eeq
For $p\geq 2$ there is connection between solution of Eq.~\eqref{eq7}
and system~\eqref{eq16} -- \eqref{eq18}.
In particular there are two lemmas.

{\bf Lemma 4.}\\
Let $u^*$ satisfy Eq.~\eqref{eq7}
and operator $\opL_1$ has left inverse. Then elements $u_i^* = \prod_{j=1}^i\opA_{2j} u^*,\, i=p,p-1, \dots, 1$
satisfy Eqs.~\eqref{eq16}, \eqref{eq17},  and  $u^*$  satisfy Eq.~\eqref{eq18}.\\

{\bf Lemma 5.}\\
Let elements $u_p^*, u_{p-1}^*, \dots, u_1^*, u^*$
satisfy Eqs~\eqref{eq16}, \eqref{eq17}, \eqref{eq18}
and operator $\opL_1$ has right inverse.
Then element $u^*$ determined from  Eq.~\eqref{eq18}
of  spit system~\eqref{eq16}, \eqref{eq17}, \eqref{eq18} is the solution 
to Eq.~\eqref{eq7}.

Proof of Lemma 4 and Lemma 5 for any natural number $p$  can be reduced to
employment  of  the skeleton chain via operators $\{ \opA_j\}_{j=1}^{2p}$
and repeats stages of the proofs of  Lemmas 2 and 3 for the case of  $p=1.$

Based on Lemmas 1--5 the following main result can be formulated. 

{\bf Main Theorem}.
{\it Let noninversible bounded operator $\opB$
has skeleton chain $$\{ B_1, \dots, B_p \},$$ operator $\opB_p\opL - \opL_1$
with definition domain in $E_p$ has bounded inverse. Let operator $\opL_1$ be defined on 
domains from   $E_i, \, i=1, \dots, p$ and $E.$ Let operator  $\opL_1$ has
inverse bounded operator. Then system \eqref{eq16}, \eqref{eq17}, \eqref{eq18}
enjoy unique solution
$\{u_p^*, u_{p-1}^*, \dots, u_1^*, u^*\},$
defined as follows
$$u_p^* = (\opB_p\opL - \opL_1)^{-1} \prod_{j=1}^p \opA_{2j} f, $$
$$u_i^* = \opL_1^{-1} \{- \prod_{j=1}^i \opA_{2j} f + \opA_{2j+1} \opL u_{i+1}^*   \}, i=p-1, \dots, 1, $$
$$u^* = \opL_1^{-1} \{ -f+ \opA_1 \opL u_1^* \}.  $$
Moreover, element  $u^*$ satisfy Eq.~ \eqref{eq7}
and $u_i^* = \prod_{j=1}^i \opA_{2j} u^*, i=1, \dots, p.$}

By setting the initial-boundary conditions to ensure the reversibility of the operators
$ \opL_1 $ and $\opB_p \opL - \opL_1$ with specific differential operators
$ \opL $ and $ \opL_1 $ and using the Main Theorem the existence and uniqueness theorems
can be derived. Moreover, the formula obtained in Theorem can effectively
build the desired classical solution of Eq.~\eqref {eq1}
with sufficient smoothness of $f: \Omega \subset {\mathbb R}^{m+1} \rightarrow E$
  and the coefficients of the differential operators $\opL$ and $\opL_1.$
Such applications of the theory are discussed further in section 3.

\section{The existence and methods of constructing\\ solutions of nonclassic BVP  with\\ partial derivatives}

\subsection{System of Eqs.~\eqref{eq19}}
Consider the system
\beq
\opB\hspace*{-1.5mm}\displaystyle{\sum\limits_{k_1+k_2 \leq n} a_{k_1 k_2} \frac{\partial^{k_1+k_2}u(x,t) }
{\partial t^{k_1} \partial x^{k_2}} }= \displaystyle{\sum\limits_{k_1+k_2 \leq m} c_{k_1 k_2} 
  \frac{\partial^{k_1+k_2}u(x,t) }
{\partial t^{k_1} \partial x^{k_2}} + f(x,t)}.
\label{eq19}
\eeq
Here $m<n,$ $\opB$ is constant  $N\times N$ matrix, $\det B=0,$
$a_{k_1 k_2},$ $c_{k_1 k_2}$ are numbers, $a_{n0 \neq 0,}$  $a_{0 n} = 0, $
$c_{0  m} \neq 0,$ $c_{m 0} =0.$  Vector-functions
$u(x,t) = (u_1(x,t), \dots, u_N(x,t))^{T},$ $f(x,t) = (f_1(x,t), \dots, f_N(x,t))^{T}$
are supposed to be defined and analytical for $-\infty < x,t < \infty.$
 
Let  $\text{rank}\, \opB =r < N.$ Then based on  \cite{gantmaher}
$\opB = \opA_1 \opA_2,$  where $\opA_1 $ is $N \times r$ matrix, 
$\opA_2 $ is $r \times N$ matrix.
Let us introduce the $r \times r$ matrix $\opB_1 = \opA_2 \opA_1$ 
and assume  $\det \opB_1 \neq 0.$
Then taking into account Lemma 3, the solution to system of Eqs~\eqref{eq19} 
can be reduced  to the successive solution of Eqs~\eqref{eq8}, \eqref{eq9},
which are as follows in this case

\beq
\opB_1 \sum\limits_{k_1+k_2 \leq n} a_{k_1 k_2} \frac{\partial^{k_1+k_2}u_1(x,t) }
{\partial t^{k_1} \partial x^{k_2}} = \sum\limits_{k_1+k_2 \leq m} c_{k_1 k_2} 
  \frac{\partial^{k_1+k_2}u_1(x,t) }
{\partial t^{k_1} \partial x^{k_2}} + \opA_2 f(x,t),
\label{eq20}
\eeq

\beq
 \sum\limits_{k_1+k_2 \leq m} c_{k_1 k_2} \frac{\partial^{k_1+k_2}u(x,t) }
{\partial t^{k_1} \partial x^{k_2}} = - f(x,t) + \opA_1 \sum\limits_{k_1+k_2 \leq n} 
a_{k_1 k_2} 
  \frac{\partial^{k_1+k_2}u_1(x,t) }
{\partial t^{k_1} \partial x^{k_2}},
\label{eq21}
\eeq
where $\det \opB_1 \neq 0,$ $u_1(x,t) = (u_{11}(x,t), \dots , u_{1r}(x,t))^{T}, r<N,$ $u_1=A_2u.$
Since by hypothesis $a_{n0} \neq 0, c_{0m}\neq 0,$  then for system~\eqref{eq19}
one may introduce the initial conditions
\beq
\frac{\partial^i u(x,t)}{\partial x^i} \biggl|_{x=0}, \, i=0,1, \dots, m-1,
\label{eq22}
\eeq

\beq
\opA_2 \frac{\partial^i u(x,t)}{\partial t^i} \biggl|_{t=0} = 0, \, i=0, 1, \dots, n-1.
\label{eq23}
\eeq

 Vector-function $u_1(x,t)$
based on Kovalevskaya theorem can be defined as unique solution to system~\eqref{eq20} 
with initial conditions
$$\frac{\partial^i u_1(x,t) }{\partial t^i} \biggl |_{t=0}=0, \, i=0,1, \dots, n-1.$$
By vector $u_1(x,t)$ substitution into the right hand side of system \eqref{eq21},
the desired vector  $u(x,t)$ can be found as unique solution to the Cauchy
problem~\eqref{eq21} -- \eqref{eq22}.

\subsection{System of Eqs.~\eqref{eq24}}

Consider system
\beq
\opB \frac{\partial^{n}u(x,t) }
{\partial x^{n} } = \biggl ( \frac{\partial}{\partial t} - a^2 \frac{\partial^2}{\partial x^2 } \biggr ) u(x,t) +f(x,t), \, n \geq 3.
\label{eq24}
\eeq
As in system~\eqref{eq19}, $\opB $ is singular $N \times N$ matrix
with $ \text{rank}\, \opB =r <N,$
$\opB = \opA_1 \opA_2,$
$\opB_1 = \opA_2 \opA_1,$ $\det B_1 \neq 0.$
 Let $f(x,t) = (f_1(x,t), \dots , f_N(x,t))^{T}$ be vector-function defined  
for $0\leq x \leq 1,$ $0< t < \infty,$  continuous with respect to  $x$ and analytical by
 $t,$  $u=(u_1,\dots , u_N)^{T}.$

The objective is to construct solution of system of Eqs~\eqref{eq24} in
$\Omega = \{ 0\leq x \leq 1, 0< t < \infty \}.$
Based on Lemma 4  and Main Theorem introduce system of two equations ($u_1=\opA_2 u$)
\beq
\opB_1 \frac{\partial^{n}u_1(x,t) }
{\partial x^{n} } = \biggl ( \frac{\partial}{\partial t} - a^2 \frac{\partial^2}{\partial x^2 } \biggr ) u_1(x,t) + \opA_2 f(x,t).
\label{eq25}
\eeq

\beq
\biggl ( \frac{\partial}{\partial t} - a^2 \frac{\partial^2}{\partial x^2 } \biggr ) u(x,t)  = - f(x,t)
+\opA_1 \frac{\partial^n u_1(x,t) }{\partial x^n}.
\label{eq26}
\eeq
with initial-boundary conditions
\beq
\frac{\partial^i u_1(x,t)}{\partial x^i} \biggl|_{x=0} = 0, i=0,1,\dots , n-1
\label{eq27}
\eeq
\beq
u(x,t)|_{t=0} = 0
\label{eq28}
\eeq
\beq
u(x,t)|_{x=0} = 0, \, u(x,t)|_{x=1}= 0.
\label{eq29}
\eeq
Since $\det B_1 \neq 0$ then vector-function $u_1(x,t)$  based
on Kovalevskaya theorem can be defined as unique solution of 
Cauchy problem~\eqref{eq25} --\eqref{eq27}.
Substitute $u_1(x,t)$ in the right hand side of  \eqref{eq26},
 the unique solution of the first boundary value problem
\eqref{eq26}, \eqref{eq28} \eqref{eq29} is constructed for heat equation
using known formula (here readers may refer to p.~215 in \cite{tikhon}), using the source function.
 Constructed solution $u(x,t)$ will be classic unique solution 
of system~\eqref{eq24} in domain $\Omega = \{ 0\leq x \leq 1, \, 0<t < \infty  \}.$
This solution
 satisfy initial conditions
$$\opA_2 \frac{\partial^i u(x,t)}{\partial x^i} \biggl|_{x=0} = 0, \, i=1,\dots, n-1$$
 and conditions~\eqref{eq28} -- \eqref{eq29}.

\subsection{Integro-differential equation of the first kind}

Consider integro-differential equation of the first kind
 \beq
\int_a^b K(x,s) \frac{\partial^3 u(s,y,t)}{\partial t^3}\, ds = \frac{\partial^2 u(x,y,t)}{\partial y^2}
+a u(x,y,t) +f(x,y,t),
\label{eq30}
\eeq
where $x\in [a,b],$ $s \in [a,b],$ $y\in [0,1],$   $0\leq t < \infty,$ $K(x,s) = \sum_{j=1}^n a_j(x) b_j(s)$ be continuous kernel, 
$f(x,y,t)$ is continuous function,  $ a\neq (\pi n)^2.$
Here operator $\opB := \int_a^b K(x,s) [\cdot ] \, ds, $ $\opL = \frac{\partial^3}{\partial t^3},$
$\opL_1 = \frac{\partial^2}{\partial y^2 }+a.$
Let $$\det \left[\int_a^b a_i(s) b_j(s)\, ds\right]_{i,j=1}^n \neq 0.$$
Then using {\it Definition 1}, matrix 
$\opB_1 = ||\int_a^b a_j(s) b_i(s)\, ds||_{i,j=1}^n$ be the regular
skeleton chain of  integral operator $\opB$ of length 1.
Then using Lemma 3 the system~\eqref{eq25}--\eqref{eq26}  can be introduced in the following form
\beq
\opB_1 \frac{\partial^3 u_1(y,t)}{\partial t^3} = \left (\frac{\partial^2}{\partial y^2} + a\right )u_1(y,t) +\beta(y,t),
\label{eq31}
\eeq
 where $$\beta(y,t) = \left ( \int_a^b f(x,y,t)b_1(x)\, dx, \dots ,  \int_a^b f(x,y,t)b_n(x)\, dx,   \right )^{T},$$
$$u_1(y,t) = (u_{11}(y,t), \dots , u_{1n}(y,t))^{T},$$
$$u_{1j}(y,t) = \int_a^b b_j(s) u(s,y,t) \, ds, \, j=1,2, \dots, n,\,
\det \opB_1 \neq 0.$$
Then Eq. \eqref{eq30} has the next form
\beq
\frac{\partial^2 u(x,y,t)}{\partial y^2} + a u(x,y,t) = -f(x,y,t) + \sum_{j=1}^n a_j(x) 
\frac{\partial^3}{\partial t^3}
u_{1j}(y,t).
\label{eq32}
\eeq
Following the Main Theorem  the initial conditions can be defined as follows
\beq
\frac{\partial^i u_1(y,t)}{\partial t^i}\biggl|_{t=0} = 0, i=0, 1, 2
\label{eq33}
\eeq
and for Eq. \eqref{eq32} the boundary conditions can be defined as follows
\beq
u(x,y,t)|_{y=0} = 0,
\label{eq34}
\eeq
 \beq
u(x,y,t)|_{y=1} = 0.
\label{eq35}
\eeq
Vector function $u_1(y,t)$ can be defined as solution 
of the regular Cauchy problem \eqref{eq31}, \eqref{eq33}.
Substitute $u_1(y,t)$ into the right hand side of Eq.~\eqref{eq32}, where  $a \neq (n \pi)^2.$
Then desired solution $u(x,t)$ of Eq.~\eqref{eq30}  can be defined as unique 
solution of boundary problem~\eqref{eq32}, \eqref{eq34}, \eqref{eq35} as follows
 $$u(x,y,t) = \int\limits_0^1 G(t,t_1) \biggl\{ -f(x,y,t_1) + \sum_{j=1}^n a_j(x) u_{1j}(y,t_1)  \biggr\} \, dt_1,$$
where $G(t,t_1)$ is the Green function.
Constructed solution $u(x,y,t)$  satisfies Eq.~\eqref{eq30}
and conditions
$$u(x,y,t)|_{y=0} = 0,\, u(x,y,t)|_{y=1} = 0, $$
$$\int_a^b b_j(s) \frac{\partial^i u(s,y,t)}{\partial t^i}\, ds \bigl|_{t=0} = 0, \, j=1,2,\dots , n, \, i=0,1,2. $$

\section{Skeleton decomposition in the theory of\\ irregular ODE in Banach space}

Consider the simplest irregular ODE
\begin{equation}
\opB \frac{d u(t)}{dt} = u(t)+ f(t),
\label{eq43n}
\end{equation}
$f(t): \, [0,\infty) \rightarrow E,$ $\opB \in {\mathcal L}(E \rightarrow E).$
Let $\{\opB_1, \dots, \opB_p \}$ is skeleton chain of operator $\opB.$
Then from Main Theorem the following next results can be formulated.\\
$\,$\\
{\it Theorem 1.\\
Let $\{\opB_1, \dots, \opB_p \}$ be regular skeleton chain, function $f(t)$
$p-1$-times differentiable then Eq.~\eqref{eq43n} with initial condition
\begin{equation}
\prod_{j=1}^p A_{2j} u(t)|_{t=0} = c_0, \, c_0 \in E_p
\label{eq44n}
\end{equation}
enjoy unique classic solution $u_0(t,c_0).$ Here 
\begin{equation}
u_0(t,c_0) = -f(t) +{\mathbf A}_1 \frac{d u_1}{dt}.
\label{eq45n}
\end{equation}
}

Let us outline the scheme for construction function $u_1(t,c_0)$ in solution \eqref{eq45n} to 
problem \eqref{eq43n} -- \eqref{eq44n} as follows:

\begin{enumerate}

\item If $p=1$ then $u_1(t,c_0)$ satisfy the regular Cauchy problem
$$ \left\{
            \begin{array}{ll}
              \opB_1 \frac{d u_1}{d t} = u_1 + {\mathbf A}_2 f(t), \\
              u_1(0) = c_0.
            \end{array}
          \right. $$

\item If $p \geq 2$ then function $u_1(t,c_0)$ can be constructed by following recursion

$$ \left\{
            \begin{array}{ll}
              \opB_p \frac{d u_p}{d t} = u_p + \prod_{j=1}^p {\mathbf A}_{2j} f(t), \\
              u_p(0) = c_0.
            \end{array}
          \right. $$

\end{enumerate}
$$u_i(t,c_0) = {\mathbf A}_{2i+1} \frac{d u_{i+1}(t,c_0)}{d t}  - \prod_{j=1}^i {\mathbf A}_{2j} f(t),\, i=p-1,p-2, \dots , 1.$$
{\it Theorem 2.}\\
Let $\{\opB_1, \dots , \opB_{p-1}, 0\}$ be singular chain of lenth $p,$
$0$ is null operator acting from $E_p$  to $E_p$ be singular skeleton chain of 
length $p\geq 1.$
Then $\opB$ be nilpotent operator and homogeneous eq. $\opB \frac{du}{dt} = u$
has only trivial solution. In this case, if function $f(t)$ be $p$-times
differentiable then unique classic solution of eq. \eqref{eq43n} can be constructed
as follows
$$ u_n(t) = -f(t) + \opB \frac{d}{dt} u_{n-1}(t), \, u_0(t) = -f(t), \, n=1,2, \dots , p. $$
Here function $u_p(t)$ be unique classic solution of eq. \eqref{eq43n}.

\subsection{Examples}

In order to illustrate Theorem 1 and Theorems 2 let us consider the following examples.

{\bf Example 1.}\\
Let in Eq.~\eqref{eq43n} 
\[ \opB = \left( \begin{array}{ccc}
0 & 2 & 2 \\
0 & 1 & 1 \\
0 & 0 & 0 \end{array} \right),\] 
$f(t) = (f_1(t), f_2(t), f_3(t)  )^{T},$ $u(t) = (u_1(t), u_2(t), u_3(t))^{T}.$
$\opB  = {\mathbf A}_1  {\mathbf A}_2,$
where $${\mathbf A}_1 = \left( \begin{array}{c}
2  \\
1  \\
0  \end{array} \right),$$
${\mathbf A}_2 = (a,\cdot)$ (scalar product),
$a = (0,1,1)^{T}.$
Since $\opB_1 = {\mathbf A}_2 {\mathbf A}_1 = 1,$
then matrix $\opB$ generates the regular skeleton chain of lenght $p =1.$
Based on Theorem 2 the following system can be written
$$ \left\{
            \begin{array}{ll}
              \frac{dv}{dt}= v +f_2(t) +f_3(t), \\
		v(0) = 0,\\
              u(t) = -f(t) +{\mathbf A}_1 \frac{dv}{dt},
            \end{array}
          \right. $$
where $v= {\mathbf A}_2 u = u_2(t) +u_3(t).$
The initial condition $u_2(0)+u_3(0) = 0$ has defined to meet condition $v(0)=0.$
It can be assumed that $f_i(t)$ is continuous function.
Then $$v(t) = \int_0^t e^{t-s} ( f_2(s) - f_3(s) )\, ds $$
and system \eqref{eq43n} is defined by matrix $\opB$
and initial condition $u_2(0)+u_3(0)=0$
enjoy unique classic solution 
$u(t) = -f(t) +(2v(t) +2 f_2(t) +2f_3(t), v(t)+f_2(t)+f_3(t), 0)^{T}.$

{\bf Example 2.}\\
Let in Eq.~\eqref{eq43n} 
\[  \opB =\left( \begin{array}{ccc}
0 & 2 & -2 \\
0 & 1 & -1 \\
0 & 1 & -1 \end{array} \right),\, {\mathbf A}_2 = \left( \begin{array}{c}
2  \\
1  \\
1  \end{array} \right), \] ${\mathbf A}_2 = (a,\cdot),$ $a=(0,1,-1)^{T}.$
Then $\opB_1 = {\mathbf A}_2 {\mathbf A}_1 = 0,$
 i.e. singular skeleton chain length $p=1$ corresponds to matrix $\opB,$
$\opB^2$ is null matrix.
Therefore according to Theorem 2 solution to the system
\eqref{eq43n} with such matrix $\opB$ can be constructed as follows
$$u_n(t) = -f(t) +\opB \frac{d}{dt} u_{n-1}(t), n=1,2,\dots , u_0(t) = -f(t). $$
For $n=1$ we have the desired solution can be constructed as follows
$$u(t) = -f(t) +(2(f_2(t)-f_3(t))', (f_2(t) - f_3(t))',   (f_2(t) - f_3(t))' )^T. $$
 If $f(t)$ is differentiable function, then we get classic solution.
If $f_i(t)$ is piecewise absolutely  continuous function with 
discontinuity points of the 1st kind and piecewise continuous derivative
then we get solution in   the space of  distributions $K'.$
Proposed method enables generalized solutions construction as well.

\section{Generalizations and additional approaches to  formulation of BVPs
for the irregular PDEs}

Suggested approach based on skeleton chain of operator $\opB$
can be employed for nonlinear equation as well. Indeed, let us replace in Eq.~\eqref{eq1} 
$f(x)$ with nonlinear by $u$ mapping $f(Mu,x)$. It is assumed 
 $M = \prod_{j=1}^p \opA_{2j},$
 where $\{\opA_{2j} \}$ are operators generate skeleton chain of operator  $\opB.$ 
Then Eq.~\eqref{eq16} of reduced system would be
nonlinear equation in $u_p.$ But in that case $\opB_ p$ is invertible in Eq.~\eqref{eq16}.
Therefore, it is easy to formulate analogues of Lemmas 1-5 and Main Theorem for such 
nonlinear case if mapping  $f(Mu,x) $ is assumed Lipschitz in $u.$

The novel theory of non-classical boundary value problems for Eq.~\eqref{eq1}
can be developed   taking into account the Jordan structure of  operator $ \opB. $
This approach has been used to  study of boundary value problems for hyperbolic
Goursat systems with singular matrix in the main part.

Using  system \eqref{eq24} it is demonstrated below that the method of work \cite{sid_blag} 
can also be applied to Eq.~\eqref{eq1}.

Indeed, let in system \eqref{eq24} matrix $\opB$ is $N\times N$  symmetric matrix, $\text{rank}\, \opB = r,\,$ and let $\{ \varphi_i \}_1^{N-r}$ be orthonormal
basis in $\ker \opB,$ $\varphi_i = (\varphi_{i\,1}, \dots , \varphi_{i\,N})^{T},$ $u = (u_1,\dots, u_N)^{T}.$

Introduce matrix $$\Gamma= (\opB + \sum_{i=1}^{N-r} ( \cdot, \varphi_i ) \varphi_i)^{-1}.$$
Then, without loss of generality, we can seek a solution of  system \eqref{eq24}
as following sum
\beq
u(x,t) = \Gamma v(x,t) + \sum_{i=1}^{N-r} c_i(x,t) \varphi_i,
\label{eq36}
\eeq
where $v = (v_1,\dots, v_n)^{T},$
\beq
( v(x,t), \varphi_i ) = 0, \, i=1,2, \dots N-r.
\label{eq37}
\eeq

Introduce projector  ${\mathbf P} = \sum_{i=1}^{N-r} ( \cdot, \varphi_i  ) \varphi_i$.
Due to \eqref{eq37}
$({\mathbf I} - {\mathbf P}) u = \Gamma v,$ $\opB \Gamma = {\mathbf I} - {\mathbf P}.$

${\mathbf P}v=0,$ ${\mathbf P} \Gamma = \Gamma {\mathbf  P}.$
Therefore, substitute  \eqref{eq36} into system\eqref{eq24}. Then application of operator
${\mathbf I -\mathbf P},$ followed by ${\mathbf P}$ yield two equations
\beq
\frac{\partial^n v(x,t)}{\partial x^n} = \left (\frac{\partial}{\partial t} - a^2 \frac{\partial^2}{\partial x^2} \right ) \Gamma v(x,t) + ({\mathbf I - \mathbf P}) f(x,t),\, n\geq 3,
\label{eq38}
\eeq

\beq
\left (\frac{\partial}{\partial t} - a^2 \frac{\partial^2}{\partial x^2} \right ) c(x,t) = \beta(x,t),  
\label{eq39}
\eeq
where $c(x,t) = (c_1(x,t), \dots , c_{N-r}(x,t))^{T},$\\
$\beta(x,t) = (( f(x,t), \varphi_1), \dots , ( f(x,t), \varphi_{N-r}))^{T}.$
Define for system~\eqref{eq24} boundary conditions
\beq
({\mathbf I - \mathbf P}) \frac{\partial^i u(x,t)}{\partial x^i}\bigr|_{x=0}=0, \, i=0,1,\dots n-1,
\label{eq40}
\eeq
\beq
{\mathbf P} u(x,t)|_{t=0} = 0, {\mathbf P }u(x,t)|_{x=0} = 0, {\mathbf P} u(x,t)|_{x=1} = 0.
\label{eq41}
\eeq
System~\eqref{eq38} is Kovalevskaya type.
Conditions \eqref{eq40} generate the following Cauchy conditions 
\beq
 \frac{\partial^i v(x,t)}{\partial x^i}\bigl |_{x=0} = 0, \, i=0,1, \dots, n-1
\label{eq42}
\eeq
for determination vector-function $v(x,t)$ from system~\eqref{eq38}.

Conditions \eqref{eq41} generate for parabolic Eqs.~\eqref{eq39}
the known boundary conditions 
\beq
c_i(x,0)=0, c_i(0,t)=0, c_i(1,t)=0, \, i=1,2, \dots, N-r.
\label{eq43}
\eeq

If vector-function  $f(x,t)$ is continuous on $t$ and analytical on  $x,$
then  $v(x,t)$ can be constructed as classic solution to Cauchy problem
\eqref{eq38}--\eqref{eq42} based on Kovalevskaya theorem. 
Vector-function $c(x,t)$ can be defined by solution of parabolic Eq.~\eqref{eq39} 
with conditions~\eqref{eq43} in closed form using source function, see p.~215 in \cite{tikhon}.

Other complicated boundary value problems for PDE systems can be raises and resolved
using  the generalized Jordan chains of linear operators.
Recent results in the theory of generalized Jordan chains of linear operators are presented 
in \cite{loginov}.

\section*{Conclusion}

This paper reports on the novel method of skeletal chains  initiated in   \cite{arxiv} for the linear operators in order to  produce new non-classical boundary value problems for systems of differential and integral-differential equations with partial derivatives arising in the modern mathematical modeling.

\end{document}